\documentclass[12pt]{article}
\usepackage{geometry}                
\usepackage{graphicx}
\usepackage{amssymb}
\usepackage{epstopdf}
\usepackage{graphicx}
\usepackage{amssymb}
\usepackage{epstopdf}
\usepackage{amsmath}
\usepackage{amsthm}
\usepackage{amsfonts}
\usepackage{epstopdf}

\newtheorem{theorem}{Theorem}
\newtheorem{lemma}[theorem]{Lemma}

\newtheorem{proposition}[theorem]{Proposition}
\theoremstyle{definition}

\theoremstyle{remark}

\numberwithin{equation}{section}

\newcommand\RR{{\mathbb R}}

\newcommand\EE{{\mathbb E}}
\newcommand\PP{{\mathbb P}}

\def\W2{W^{1,2}({\cal O}(M))}

\title{Hard edge tail asymptotics} 

\author{{Jos\'e A. Ram\'irez}\footnote{Department of Mathematics,  Universidad de Costa Rica. e-mail:
{\tt{alexander.ramirez$\_$g@ucr.ac.cr.}}},  { Brian Rider}\footnote{Department of Mathematics, University of Colorado Boulder. e-mail: {\tt{brian.rider@colorado.edu.}}}, { Ofer Zeitouni}\footnote{ 
Department of Mathematics, University Minnesota; 
Department of Mathematics, Weizmann Institute. e-mail:
{\tt{zeitouni@math.umn.edu.}} }
}
\date{} 

\begin{document}
\maketitle

\begin{abstract} 
Let $\Lambda$ be the limiting smallest eigenvalue in the general $(\beta, a)$-Laguerre
ensemble of random matrix theory.  That is,  $\Lambda$ is the $n \uparrow \infty$ distributional limit of  
the (scaled) minimal point drawn from the  density proportional to 
$\prod_{1 \le i < j \le n} | \lambda_i- \lambda_j|^{\beta} \, $
$\prod_{i=1}^n  \lambda_{i}^{\frac{\beta}{2}(a+1) -1}  e^{-\frac{\beta}{2} \lambda_i} $  on $(\RR_+)^n$.  Here $\beta>0$, $a >-1$;
for $\beta=1,2,4$ and integer $a$, this object governs the singular values of certain rank $n$  Gaussian matrices. We prove that 
$$
  P(\Lambda > \lambda) =    e^{- \frac{\beta}{2} \lambda + 2 
\gamma \sqrt{\lambda}} \lambda^{ - \frac{\gamma(\gamma+1-\beta/2)}{2\beta} 
   } \, {\mathfrak{e}} (\beta, a)  (1+o(1)) 
$$
as $\lambda \uparrow \infty$ in which $\gamma = \frac{\beta}{2} (a+1)-1$ and ${\mathfrak{e}}(\beta, a)>0$ is a constant (which we do not determine).  This estimate complements/extends various results previously available for special values of $\beta$ and $a$.
\end{abstract}

\section{Introduction}

The shape of the distribution of the smallest singular value of a ``typical" matrix is a deeply studied question.  An overview of the varying motivations for this problem may be found in \cite{RV}.  In the case of Gaussian matrices, many exact formulas are available both at finite dimension and asymptotically \cite{Edelman, TW}.  Only quite recently has it been shown that the asymptotic laws are universal beyond the Gaussian case (in the sense of being insensitive to the statistics of the matrix entries), see \cite{TV}.

\medskip

Here we consider the ``general beta" analogues of the classical Gaussian ensembles.   These are defined by placing a measure on $n$ nonnegative real points $\lambda_1, \lambda_2, \dots, \lambda_n$ with density function
(a normalization constant times)
\begin{equation}\label{density}
  \prod_{1 \le i < j \le n} | \lambda_i- \lambda_j |^{\beta} \,  \prod_{i=1}^n  \lambda_{i}^{\frac{\beta}{2}(a+1) -1} 
  e^{-\frac{\beta}{2} \lambda_i} .
\end{equation}
When $\beta=1,2,4$ and $a = 0, 1, 2, \dots$  this is the joint 
square-singular value law of an $n \times (n+a)$ real, complex, or quaternion Gaussian matrix.  It is however a sensible law for any $\beta > 0$ and $a > -1$, and, what is more, still a joint square-singular value law for a certain random bi-diagonal matrix ensemble \cite{DE}.  Further, the least order statistic $\lambda_{min}$ satisfies a limit law:  as $n \uparrow \infty$,  $n^2 \lambda_{min}$ converges (in distribution) to a well defined random variable, denoted here by $\Lambda$ ($=\Lambda(\beta, a)$).
There are several proofs of this for special values of $\beta$ and $a$; \cite{RR} contains a proof (making use of the
bi-diagonal representation of \cite{DE} and substantiating a conjecture of \cite{ES}) valid for all values of those parameters.

\medskip

Our starting point is a relation between the law of $\Lambda$ and the explosion/non-explosion of the diffusion process: with $b$ a Brownian motion,
\begin{equation}\label{maindiff}
d x(t) = db(t) + \left( \frac{\beta}{4} (a + \frac{1}{2} ) -  \frac{\beta}{2} \sqrt{\lambda}  e^{-\beta t /8} 
                                                    \cosh x(t) \right) dt.
\end{equation}
In particular, a corrected version of Theorem 2 of \cite{RR} (see also the derivation leading to (\ref{varchange}) below) implies that
\begin{equation}
\label{passage}
  P(\Lambda > \lambda) = \PP_{\infty, 0} ( t \mapsto x(t) \mbox{ never hits } - \infty ). 
\end{equation}
Here $\PP_{c, s}$ indicates the law on paths induced by $x$, begun from position $c$ at time $s$.  Our main result reads:

\begin{theorem}
\label{main}
 Let ${\mathfrak p}_{\lambda} = {\mathfrak p}_{\lambda, \beta, a}$ denote the right hand side of (\ref{passage}).  For large values of $\lambda$ it holds
\begin{equation}\label{result}
   {\mathfrak p}_{\lambda} = e^{- \frac{\beta}{2} \lambda + 2 \gamma 
\sqrt{\lambda}} \lambda^{ - \frac{\gamma(\gamma+1-\beta/2)}{2\beta} 
} \, {\mathfrak e}(\beta, a)  (1+o(1)).  
\end{equation}
Here $\gamma = \frac{\beta}{2}(a+1)-1$ and ${\mathfrak e}(\beta, a) > 0$ is an undetermined constant.
\end{theorem}
 
There has already been a great deal of work in this direction, though focussed on dealing directly with the statistics (\ref{density}) rather than our passage time description (\ref{passage}).  The fundamental treatment of Tracy-Widom \cite{TW}
for $\beta=2$ produced the correct $\lambda \rightarrow \infty$ asymptotics of $\mathfrak{p}_{\lambda, 2, a}$ up to a multiplicative constant 
and provided a conjecture for that constant, ${\mathfrak{e}}(2, a)$.  This  has recently been verified by Ehrhardt \cite{Torsten},  for $|a| <1$, by operator theoretic techniques, and for all $a > -1$ by Deift-Krasovsky-Vasilevska \cite{DKV}  using Riemann Hilbert Problem machinery.  
A non-rigorous  argument in \cite{CM} predicted all factors in the asymptotics save the constant  for all $(\beta, a)$.  Making use of integral identities available at special values of $\beta$ and integer $a$, Forrester has a sound conjecture for the value of the general constant ${\mathfrak e}(\beta, a)$, see \cite{For}.  It appears possible that a type of analytic continuation argument could extend the result of \cite{For} to all $\beta$ (though for still integer $a$).  The method employed here leaves 
${\mathfrak e}(\beta, a)$ in opaque form, as a somewhat involved expectation over diffusion paths; an explicit determination of this object for all $\beta$ and $a$ remains an open problem.

\medskip

In many ways, the chief insight of this paper is to cast the diffusion (\ref{maindiff}), which encodes the desired probability distribution, in the present form. (The process which appears in \cite{RR} is related by a change of variables.)  In fact, $t \mapsto x(t)$ is remarkably similar to the process studied by Valk\'o-Vir\'ag in estimating the probability of large gaps in the general beta ``bulk" \cite{VV}.  They showed that the probability of a gap being larger than $\lambda$ is equal to the 
non-explosion, again to $-\infty$,  of 
$$
  d z(t) = db(t)  + \left( \frac{1}{2} \tanh z(t) -  \frac{ \beta}{8}  \lambda  e^{-\beta t/4} \cosh z(t)  \right) dt, 
$$
begun again at $+\infty$.  It is no surprise then that their basic argument, which involves estimating the Cameron-Martin-Girsanov factor produced by a well-chosen change of measure, may be followed in this case.  

\medskip

The proof of Theorem 1 occupies sections 3 and 4; section 2 gives a self-contained explanation of the identity (\ref{passage}).

\section{Passage time description for $\Lambda$}

Without pointing the reader to \cite{RR} and the subsequent erratum, it is easy enough to give a brief derivation of the relevance of the diffusion (\ref{maindiff}) to the distribution function $P(\Lambda > \lambda)$.  The main result of
\cite{RR} shows that $\Lambda^{-1}$ is the maximal eigenvalue of the almost surely trace class integral operator
\begin{equation}
\label{L_integral}
      L_{\beta, a} \psi(t) :=  \int_0^{\infty} \left(    \int_0^{t \wedge s} e^{au + \frac{2}{\sqrt{\beta}} b(u)} \, du \right)    \psi(s)  
        e^{-(a+1)s - \frac{2}{\sqrt{\beta}} b(s)} \, ds, 
\end{equation}
acting on $L^2[ \RR_+, \mu]$, $\mu(dt) = e^{-(a+1)t - \frac{2}{\sqrt{\beta}} b(t)} \, dt$. Here $t \mapsto b(t)$ is a standard Brownian motion.

\medskip

Any nonnegative $L^2$ solution of $\psi(t) = \lambda L_{\beta, a} \psi(t)$ 
satisfies $\psi(0)=0$ and $ \psi'(t) :=d\psi(t)/dt\ge 0$
for all $t > 0$, as can be seen by taking derivatives of both sides of the eigenvalue equation:
\begin{equation*}\label{Neumann}
  \psi'(t) = \lambda  e^{at + \frac{2}{\sqrt{\beta}} b(t)} \int_t^{\infty}   \psi(s)  
        e^{-(a+1)s - \frac{2}{\sqrt{\beta}} b(s)} \, ds.
\end{equation*}
This converts to a differential system:  
\begin{equation}\label{system} 
  d \psi'(t) =  \frac{2}{\sqrt{\beta}} \psi'(t) db(t) + [ (a+\frac{2}{\beta}) \psi'(t) - \lambda e^{-t} \psi(t) ] dt, \   \    d\psi(t) = \psi'(t) dt
\end{equation}
which can be used to test whether a fixed $\lambda$ is at or below an eigenvalue.  Specifically, $\lambda$ is strictly below
the groundstate eigenvalue $\Lambda$ if the solution to (\ref{system}) begun at $\psi(0)=0$ (and $\psi'(0) = 1$ say) satisfies 
$\psi(t) > 0, \psi'(t) > 0$ for all time (note that solutions of (\ref{system}) are decreasing in $\lambda$).   It is now the standard trick to translate this condition onto the diffusion $q(t) := \psi'(t)/ \psi(t)$ which
solves
$$
   dq(t) = \frac{2}{\sqrt{\beta}} q(t) db(t)  + [ (a+ \frac{2}{\beta}) q(t) - q^2(t)   - \lambda e^{-t} ] dt,
$$
started from $+\infty$ at time $t=0$.   In particular, if $\tau_c$ is the passage time of $q$ to a level $c$, the event
$\{ \Lambda > \lambda \}$ coincides with $\{ \tau_0 = \infty\}$. Now the change of variables,
\begin{equation}\label{varchange}
  x(t) :=  \log ( q(\beta t/4) ) + \beta t/8 - \log \lambda/2,  
\end{equation}
explains the identity (\ref{passage}).

\medskip

As a bit of amplification, we remark that  for $q = q(\cdot; a, \beta, \lambda)$ with $a \ge 0$,  
\begin{equation}
\label{boundaries}
\PP( \tau_{-\infty}(q) < \infty | \tau_0(q) < \infty) =1.
\end{equation}
So, at least for $a \ge 0$, one can replace the condition of $q$ never vanishing with the (more familiar) condition that $q$ never explodes to $-\infty$. 
Furthermore, a change of variables similar to (\ref{varchange}) 
shows that the event that $q(\cdot; a , \beta, \lambda)$ started from $0$ never hits $-\infty$ is the same as the event that $q(\cdot; -a-1, \beta, \lambda)$, started from $+\infty$ never hits $0$.  In particular, defining instead $x(t) :=  - \log ( -q(\beta t/4) ) + \beta t/8 - \log \lambda/2$ (keeping in mid that  for all $t > \tau_0$, $q(t) < 0$), $x(t)$ will solve (\ref{maindiff}) with $a$ replaced by $-a-1$.  One concludes 
that 
$$ \lim_{a \downarrow -1} P(\Lambda > \lambda) = 
\lim_{a \downarrow  -1} \PP_{\infty} (\tau_{0}(q(\cdot; a, \beta, \lambda)) = \infty) = 
\lim_{a \uparrow 0} \PP_{0} (\tau_{-\infty}(q(\cdot; a, \beta, \lambda)) = \infty) 
 = 0$$ for any $\lambda > 0$ (by say monotone convergence), as would have been guessed ahead of time.

\medskip

To prove (\ref{boundaries}), on the event $\{ \tau_0< \infty\}$ introduce the simpler change of variables $u(t) = \log(-q(t +\tau_0))$. This process
satisfies 
$$
   du(t) = \frac{2}{\sqrt{\beta}} db(t)  + [ a+  e^{u(t)}  + \lambda e^{-\tau_0} e^{-t} e^{-u(t)} ] dt,  \  \  u(0+) \in (-\infty, \infty),
$$
to which we compare the homogeneous  process defined by
$$
   dv(t) = \frac{2}{\sqrt{\beta}} db(t)  + [ a+  e^{v(t)}]  dt,  \  \   v(0) = u(0+) \in (-\infty, \infty).
$$
As $u(t) > v(t)$, $q$ explodes to $-\infty$ in finite time if $v$ explodes to $+\infty$ in finite time (we continue to work 
on the event $\{ \tau_0(q) < \infty \}$). 

\medskip

Now apply Feller's test, in the form given by Proposition 5.32 (part (ii)) 
of \cite{KS}. In particular, bring in the Lyapunov function 
$$
    m(x) = \int_0^{x} s(y) \int_0^y \frac{1}{s(z)} dz dy \   \mbox{ where }  \  s(x) = \exp \left( - \frac{\beta}{2} \int_0^x ( a + e^z) dz \right)  
$$
($s(x)$ is the derivative of the scale function for $v$). Since 
$$
   \lim_{x \rightarrow \infty}  m(x) < \infty,  \mbox{ while, if } a \ge 0, \lim_{x \rightarrow -\infty} m(x) = -\infty, 
$$
the cited form of Feller's test implies that $S = \int \{ t : v(t) \notin (-\infty, \infty)\}$ is finite with probability one. However, it is impossible that
$v(t)$ ever hits $-\infty$ (it is easily bounded below by a Brownian motion with constant drift $a$).  This completes the proof.

\section{Change of measure}

Hereafter it is convenient to put the time index in subscripts, {\em{i.e.}}, $x(t)$ becomes $x_t$ and the like.
To begin, introduce the notation
$$
  {\mathfrak p}_{\lambda}(c) = \PP_{c}( x_t \mbox{ never explodes} ).
$$
Then, by the strong Markov property,
\begin{equation}
\label{s1}
 {\mathfrak p}_{\lambda} =  {\mathfrak p}_{\lambda}(\infty) = \EE_{\infty} \Bigl[   {\mathfrak p}_1(x_T),  x_t  >  - \infty \mbox{ for } t \in (0, T] \Bigr],
\end{equation}
upon choosing 
\begin{equation}
\label{timechoice}
T =  \frac{4}{\beta} \log \lambda.
\end{equation}  
The change of measure is now enacted on the expectation (\ref{s1}).

\begin{proposition}  Let $h(t,x)$ be $C^1$ in both variables and bounded for $t \le T$. Then, the law on paths up to time $T$ 
induced by
$$
   dy_t = db_t + \left(  h(t, y_t) -  \frac{\beta}{2} \sqrt{\lambda}  \,  e^{-\beta t/ 8}  
                                                    \sinh(y_t) \right) dt,  \  \  y_0 = \infty
$$ 
is absolutely continuous with respect to that of $t \mapsto x_t$, $x_0= \infty$, subject to $x_t > -\infty, 0 \le t \le T$. 
Moreover,
\begin{equation} 
\label{s2}
{\mathfrak p}_{\lambda} = \EE_{\infty} [ {\mathfrak p}_1(y_T) \, {R_T(y_\cdot)} ],
\end{equation}
in which, for $s\leq T$,
$$
  \log R_{s}(y_.) =   \int_0^s (f(t, y_t) - g(t, y_t) ) d y_t -  
  \frac{1}{2} \int_0^s (f^2(t, y_t) - g^2(t, y_t)) dt,   
$$
$ f(t,y) = \frac{\beta}{4} (a + \frac{1}{2} ) -  \frac{\beta}{2} \sqrt{\lambda}  \,  e^{-\beta t/8}  
                                                    \cosh y $ and $g(t,y)=  h(t, y) -  \frac{\beta}{2} \sqrt{\lambda}  \,  e^{-\beta t/ 8}  
                                                    \sinh y $.
\end{proposition}

\medskip

This is just the formula of Cameron-Martin-Girsanov, applied to the particular case of a diffusion with explosion for which it is important to point out that the test function
${\mathfrak p}_1(x_T)$ in question vanishes when  $T$ 
is larger than the explosion time. One also notes that the general form of the $y$-drift,  effectively a bounded function minus
$ \sinh y$, allows $y_t$ to be started at $+\infty$ and prevents $y_t$ from exploding on $[0,T]$.  To then carry out the standard proof of Cameron-Martin-Girsanov in the present context, it must be checked that 
$R_t^{-1}(x)$, which is a local martingale by construction, is actually a martingale. But again by the general form of the $y$-drift, 
both $f-g$ and $f^2 - g^2$ are bounded when the path is bounded below, keeping $R_t^{-1}(x)$ bounded prior to the explosion time of $x_t$. Plainly, 
$R_t^{-1} (x) = 0$ at and after the explosion time.

\medskip

The first ingredient of the proof of Theorem \ref{main} is the following. Throughout the below, $[x]^{-}$ denotes the negative part of $x \in \RR$.  Also recall that
$\gamma = \frac{\beta}{2} (a+1)-1$.

\medskip

\begin{lemma} 
\label{lem1}
There exists a choice of $h$ in Proposition 2 so that,
for appropriate $\nu, \phi$ satisfying  $| \nu(t, y) | \le \kappa_1 + \kappa_2 [y]^{-}$ for all $t \ge 0$  and constant $\kappa_1, \kappa_2$, 
and $| \phi(t, y)| \le  \hat{\phi}(t) $  with $\int_0^{\infty} \hat{\phi}(t) dt < \infty$, it holds that
\begin{align}
\label{r-form}
  \log R_T(y_\cdot) 
  & =   - \frac{\beta}{2} \lambda + 2 \gamma \sqrt{\lambda}   -  
                               \left(  \frac{ \gamma ( \gamma +1 - \beta/2)}{2\beta}  \right) \log \lambda      \\
 &     +  \frac{\beta}{2} e^{-y_T} + \nu(T, y_T)  +  \int_0^T \phi(T-t, y_t) dt. \nonumber
\end{align}
\end{lemma}

\medskip
 
Once $h$ is in hand, the lemma is readily verified.  In particular,
 \begin{equation}\label{hform}
   h(t,y) =  \frac{\beta}{4} (a + \frac{1}{2})  +  h_1(y) +  e^{- \frac{\beta}{8} (T-t)} h_2(y)
\end{equation}
 where
 \begin{align}
 \label{h1h2}
    h_1(y)  & =   -  \frac{\gamma}{1+e^y},   \\
 h_2(y)   & =     \frac{1}{ \beta \sinh(y)} \left( ( h_1^2(y) - h_1^2(0) ) + \frac{\beta}{2} (a + \frac{1}{2}) (h_1(y) - h_1(0))
                    +  (h_1'(y) - h_1'(0)) \right). \nonumber
 \end{align}
 That both $h_1$ and $h_2$ are uniformly bounded, $h_1$ being integrable at $+\infty$ while $h_2$ is integrable at both $\pm \infty$
 figure into the bounds on $\nu$ and $\phi$ in the lemma.   
 
 \medskip
 
It is more instructive however to describe how $h$ is discovered, each step achieving successive order in $\lambda$, $\sqrt{\lambda}$, etc., and the various bounds claimed in the lemma seen along the way.

\medskip

{\em Step 1} begins by expanding out the exponential $R_T$ factor with a generic $h$:
\begin{align}
\label{firstpass}
\log R_T(y_\cdot) & =    -  \frac{\beta^2}{8} \lambda  \int_0^T e^{-\beta t/4} dt 
  + \frac{1}{2} \int_0^T h^2(t, y_t) dt - \frac{\beta^2}{32} (a + \frac{1}{2})^2 T  \\
  &   - \frac{\beta}{2} \sqrt{\lambda}  \int_0^T e^{-\beta t/8} h(t, y_t)  \sinh(y_t) dt
          + \frac{\beta^2}{8}    (a + \frac{1}{2} )  \sqrt{\lambda}  \int_0^T e^{-\beta t/8} \cosh(y_t) dt   \nonumber \\
   &   - \frac{\beta}{2} \sqrt{\lambda} \int_0^T e^{-\beta t/8} e^{-y_t} d y_t 
           - \int_0^T [ h(t, y_t) - \frac{\beta}{4} (a+\frac{1}{2}) ] dy_t.      \nonumber  
\end{align}
By the choice of $T$, the first term equals $- \frac{\beta}{2} (\lambda-1)$  which already gives the leading order and explains the particulars of the $\sinh y$ term in the choice of the $y$-drift.  The last term, coupled with the fact that $y_0=\infty$, prompts a
natural shift of $h$ by the factor $\frac{\beta}{4} (a+\frac{1}{2})$. That is, $h$ is replaced with $h +  \frac{\beta}{4} (a+ \frac{1}{2})$.

\medskip

{\em Step 2} enacts the above shift, and also introduces the obvious It\^o substitution in the second last term of (\ref{firstpass}),
$$
  \frac{\beta}{2} \sqrt{\lambda} \int_0^T e^{-\beta t/8} e^{-y_t} d y_t  = - \frac{\beta}{2} e^{-y_T}  + 
            ( \frac{\beta}{4} - \frac{\beta^2}{16} ) \sqrt{\lambda} \int_0^T e^{-\beta t/8} e^{-y_t} dt
$$
to write:
\begin{align}
\label{secondpass}
\log R_T(y_\cdot) & =   - \frac{\beta^2}{8} \lambda \int_0^T e^{-\beta t/4} dt    \\
  &   +  \frac{\beta}{2} \sqrt{\lambda} 
                  \left(
                 \frac{\gamma}{2}  \int_0^T e^{-\beta t/8} e^{-y_t} dt - 
          \int_0^T e^{-\beta t/8}  h(t, y_t) \sinh(y_t) d t \right)     \nonumber  \\
   &   + \frac{1}{2} \int_0^T h^2(t, y_t) dt + \frac{\beta}{4} (a + \frac{1}{2}) \int_0^T h(t, y_t) dt    - \int_0^T h(t, y_t) d y_t   + \frac{\beta}{2} e^{-y_T}. \nonumber       
\end{align}
This draws attention to line two of (\ref{secondpass}), which should produce the final constant times $\sqrt{\lambda}$ term. This may be achieved most easily by introducing a deterministic integrand in that line  via the substitution
\begin{equation}
\label{firsth}
      h(t,  y)  =  \frac{\gamma}{2} \frac{e^{-y} - 1}{\sinh(y)} + \bar{h}(t, y) := h_1(y) + \bar{h}(t, y),
\end{equation}
so that 
\begin{align*}
\lefteqn{ \hspace{-2cm}
  \frac{\gamma}{2}  \int_0^T e^{-\beta t/8} e^{-y_t} dt - 
          \int_0^T e^{-\beta t/8}  h(t, y_t) \sinh(y_t) d t } \\
          & =   \frac{\gamma}{2} \int_0^T e^{-\beta t/8} dt  - \int_0^T e^{-\beta t/8} \bar{h}(t, y_t) \sinh(y_t) d t.
\end{align*}
Evaluating all deterministic factors thus far, step 2 is summarized by
\begin{align}
\label{thirdpass}
\log R_T(y_.)   & =       - \frac{\beta}{2} \lambda + 2 \gamma \sqrt{\lambda}    + \frac{\beta}{2} e^{-y_T}  - ( {\beta} (a+ \frac{1}{2} ) + 2)\\
  &   - \frac{\beta}{2} \sqrt{\lambda} \int_0^T e^{-\beta t/8} \bar{h}(t, y_t) \sinh(y_t) dt \nonumber \\
  &   
   + \frac{1}{2} \int_0^T h^2(t, y_t) dt + \frac{\beta}{4} (a + \frac{1}{2}) \int_0^T h(t, y_t) dt    - \int_0^T h(t, y_t) dy_t.  \nonumber 
\end{align}
The first two terms above exhibit the proposed order $\lambda$ and order $\sqrt{\lambda}$ factors in the statement of the lemma, showing  that there was not much flexibility in the choice of the (uniformly bounded) function $h_1$ in (\ref{firsth}).

\medskip

{\em Step 3} is to pin down the $\log \lambda$ factor in the exponent (or, equivalently, the $T$ factor).  A look at line two of 
(\ref{thirdpass}) suggests a prescription for $\bar{h}$:
\begin{equation}
\label{secondh}
    \bar{h}(t, y) =   \frac{2}{\beta \sqrt{\lambda}} e^{\beta t/8}   h_2(y)   
                      =  \frac{2}{\beta} e^{ - (\beta/8) (T-t)}  \frac{h_3(y)}{\sinh(y)},
\end{equation}                      
in which $h_3$ must be chosen so that $h_2$ is bounded (and more).

\medskip

With $\eta(t) =  \frac{2}{\beta} e^{ - \beta t /8 } $, we employ It\^o's lemma once more to write the final term in (\ref{thirdpass}) as in
\begin{align}
\label{Ito2}
   \int_0^T h(t, y_t) dy_t &  =  H_1(y_T) + H_2(T,  y_T)  \\
    &     - \frac{1}{2} \int_0^{T} h_1^{\prime}(y_t) dt  + \int_0^T \eta'(T-t) [ \int_0^{y_t} h_2(z) dz] dt - \frac{1}{2} \int_0^{T}  \eta(T-t) h_2^{\prime}(y_t) dt.    
                 \nonumber
\end{align}
Here $H_1$ and $H_2$ denote the anti-derivative terms which appear: 
$H_1(y_T)$
 $= $ $\left.  \int_{0}^{y_t} h_1(y) dy \,  \right|_0^T $ and $H_2(T,y_T) $ $=$ $ \left. \eta(T-t) \int_0^{y_t} h_2( y) dy \,  \right|_0^T$.
Note that the boundary values of $H_2$ will necessitate that our choice of $h_2$, like that of $h_1$, is integrable at $+\infty$ ($=y_0$).  Now expand out the last two lines of (\ref{thirdpass}) to read:
\begin{align}
\label{finalexpand}
\lefteqn{  \int_0^T \left[ \frac{1}{2} h_1^2(y_t) + \frac{1}{2} h_1^{\prime}(y_t) + \frac{\beta}{4} (a+ \frac{1}{2}) h_1(y_t) -h_3(y_t) \right] dt 
    + H_1(y_T) + H_2(T,  y_T) }    \\
 &   + \int_0^T \left[ \frac{1}{2} \eta^2(T-t) h_2^2(y_t) \right.   \nonumber \\ 
 &  \quad \quad \quad
      +  \left. \eta(T-t) \left( h_1(y_t) h_2(y_t) 
                                        + \frac{\beta}{4}(a+ \frac{1}{2})    h_2(y_t) + \frac{\beta}{8}  [ \int_0^{y_t} h_2(y) dy]  + \frac{1}{2}  h_2'(y_t) \right)
 \right] dt. \nonumber
\end{align}
 The first term of (\ref{finalexpand}) prompts the choice of $h_3$, namely
set
$$
  h_3(y) = \frac{1}{2} h_1^2(y) +    \frac{1}{2} h_1'(y)+   \frac{\beta}{4} (a+\frac{1}{2}) h_1(y)  - \kappa,
$$
for a constant $\kappa$ which makes $h_2(y) = h_3(y)/\sinh(y)$ bounded.   We 
find that
$$
  \kappa=  \frac{1}{2} h_1^2(0) +   \frac{1}{2} h_1'(0)+ 
  \frac{\beta}{4} (a+\frac{1}{2}) h_1(0)   
       = \frac{\gamma(\gamma+1)}{8} - \frac{\beta}{8} ( a + \frac{1}{2}) \gamma
          = - \frac{ \gamma}{8} ( \gamma + 1 - \beta/2),
$$
compare (\ref{h1h2}).
In other words, with this choice the first term of (\ref{finalexpand}) equals $\frac{4}{\beta} \kappa \log \lambda$, the advertised $\log \lambda$ contribution in Theorem \ref{main}. 

\medskip

To finish the proof of the lemma we further identify 
\begin{equation}
\label{nu}
\nu(T, y) =    - (\beta (a + \frac{1}{2} ) + 2)   + H_1(y) + H_2(T, y),
\end{equation}
and 
$\phi(T-t, y)$ with the integrand of the last term in (\ref{finalexpand}).  One now checks:  $h_1$ and $h_2$ along with their derivatives are uniformly bounded over the entire real line (with constants depending on $a$ and $\beta$ of course), $h_1$ is integrable at $+\infty$, and $h_2$ is in fact integrable at both $\pm \infty$. The boundedness and integrability (at $+\infty$) of $h_1$ and $h_2$ give immediately that
$|H_1(y)| + |H_2(t, y)|$ $\le$ $ \kappa_1 + \kappa_2 y^{-}$, proving the claimed bound on $\nu$.
For $\phi$ one shows that  $|\phi(T-t, y)| \le c \eta(T-t)$,  using the additional appraisals on the derivatives, that $\int_{-\infty}^{\infty} |h_2|  < \infty$, and the simple fact $0 < \eta(T-t) \le \eta(0)$ for $t \in [0,T]$.

  \section{Constant term}
  
  The conclusion of the previous section is that
  \begin{eqnarray*}
  {\mathfrak p}_{\lambda}   =  
  e^{- \frac{\beta}{2} \lambda + 2 \gamma \sqrt{\lambda} }  
  \lambda^{ - \frac{\gamma(\gamma+1 - \beta/2)}{2 \beta}) \gamma } 
   \, \mathfrak{e}_{\lambda} 
  \end{eqnarray*}
 with
 \begin{equation}
 \label{elambda}
   \mathfrak{e}_{\lambda} 
  = \EE_{\infty} \left[ {\mathfrak p}_1(y_T) 
  e^{ \frac{\beta}{2} e^{-y_T} + \nu( y_T)  +
  \int_0^T \phi(T-t, y_t) dt}  \right],
 \end{equation}
  and $\nu$ and $\phi$ satisfying the bounds outlined in Lemma \ref{lem1}.
  It remains to show that the existence of a (non-zero) constant ${\mathfrak e}= {\mathfrak e} (a, \beta)$ such that
  $
    \lim_{\lambda \rightarrow \infty}  \mathfrak{e}_{\lambda}  = {\mathfrak e}.
  $
  This is again structurally identical to \cite{VV}.
  
The first observation is that the $\EE_{\infty}$ integration is performed over paths that are monotonically decreasing in $T$. The nicest way to see this is to replace the integration over $y_t,  0 \le t \le T$ with that over
$$
    {y}_t^T = y_{t+T}, \  \  -T \le t \le 0  
$$ 
which satisfies
$$
   d  {y}_t^T = db_t + ( h(t+T,  {y}_t^T) - \frac{\beta}{2} e^{-t} \sinh {y}_t^T ) dt,  \  \   {y}_{-T}^T = \infty.
$$   
If this family of processes is run on the same Brownian motion, $t \mapsto b_t$, it follows that $y_t^{T_1} \le  y_t^{T_2}$ for $t \ge -T_2$:
by definition $y_{-T_2}^{T_1} <   y_{-T_2}^{T_2}$ and the evolution maintains the ordering.
Denote this sequence of corresponding expectations simply by ${\bold{E}}$ and record that 
\begin{equation}
\label{newexp}
    \mathfrak{e}_{\lambda}
    =   {\bold{E}} [ {\mathfrak  p}_1(y_0^T)   
    e^{ \psi(y^T) }], \quad \psi(y^T) =   \frac{\beta}{2} e^{-y_0^T } +  \nu(T, y_0^T)  +  \int_0^T \phi(t, y_{-t}^T) dt .
\end{equation}

Next, pick a constant $h_0$ such that
$$
      \inf_{-\infty< y < \infty, -T < t < 0} h(t+T, y) > h_0
$$   
(a look at (\ref{hform}) and (\ref{h1h2}) shows this is possible), and introduce the stationary diffusion $t  \mapsto z_t$ on the negative half-line 
with generator
$$
  {\mathcal L} = \frac{1}{2} \frac{d^2}{dz^2} + f(z)  \frac{d}{dz}, \quad  f(z)  = h_0 - \frac{\beta}{2}  \sinh z,
$$ 
and reflected (downward) at the origin. In particular, for all $t \ge {-T}$, $\PP(z_t \in dz) = {\mathfrak{m}}(dz)  $ where
\begin{equation}
\label{zlaw}
     {\mathfrak{m}}(dz) =   \kappa_0 e^{2 h_0 z - \beta \cosh z} \, dz,  \quad z \in (-\infty, 0],  
\end{equation}
and $\kappa_0$ is the appropriate normalizer.  This is the well-known formula for the speed measure (see for example \cite{IM}), or one may check
that $\int_{-\infty}^0  {\mathcal L}  \phi(z)  \mathfrak{m}(dz) = 0 $ for all smooth $\phi$ satisfying $\phi'(0) = 0$. 

\medskip

Again running $z_t$ on the same Brownian motion, it holds that
$y_t^T \ge z_t > -\infty$ for all $t \in [-T, 0]$.  This is plain at the starting time, and continues by the domination (from below) of the  $y^T$-drift by that of $z$.  It follows that there exists a random variable $y_t^{\infty} > -\infty$ such that
\begin{equation}
\label{pointwise}
    \lim_{T \rightarrow \infty} y_t^{T}  = y_t^{\infty}  \mbox{ almost surely for each } t \in (-\infty, 0]. 
\end{equation}
To pass this convergence under the ${\bold{E}}$-expectation we prepare the following (and defer the proof to the end of the section).

\begin{lemma} 
\label{lem2}
The function $x \mapsto {\mathfrak p}_1(x)$ is continuous, strictly positive on $x > -\infty$ and satisfies
\begin{equation}
\label{hittingbound}
   {\mathfrak p}_1(x) \le  \kappa_3 e^{ -\frac{\beta}{4} e^{-x}}
\end{equation}
for a constant $\kappa_3$.
\end{lemma}

Courtesy (\ref{pointwise}) and the first statement in Lemma  \ref{lem2} we have that
\begin{eqnarray}
\label{pathwise}
    \lim_{T \rightarrow \infty} {\mathfrak p}_1(y_0^T) e^{\psi(y^T)}  & =  & {\mathfrak p}_1(y_0^{\infty})  e^{ \frac{\beta}{2} e^{-y_0^{\infty} } +  
   \nu(\infty, y_0^{\infty})  +  \int_0^{\infty} \phi(t, y_{-t}^{\infty}) dt} \\
     &:=& {\mathfrak p}_1(y_0^{\infty}) e^{\psi_{\infty} (y^{\infty})}, \nonumber
\end{eqnarray} 
using continuity (for the first three factors) and dominated convergence (for the last factor).  The evaluation $\nu(t, y) |_{t=\infty}$
simply has the effect of setting of one $H_2$-terms of which $\nu$ is comprised to zero, recall $(\ref{nu})$.

\medskip

Next, by the estimates on $\nu$, $\phi$ from  Lemma \ref{lem1} and (\ref{hittingbound}), there are the bounds
\begin{eqnarray}
\label{updown}
    \kappa_4^{-1} {\mathfrak p}_1(y_0^{T}) e^{ - \kappa_5 [y_0^{T}]^{-} }   \le    
    {\mathfrak p}_1(y_0^{T}) e^{\psi (y^{T})}
     \le  
    \kappa_4  e^{  \kappa_5   [y_0^{T}]^{-}   + ( \frac{\beta}{2}  -  \frac{\beta}{4} ) e^{-[y_0^{T}]^{-} }},
\end{eqnarray}
with positive constants $\kappa_4, \kappa_5$.  

\medskip

Note that both bounds in (\ref{updown}) depend only on the marginal of the process at time $0$, 
and denote the left and right hand sides by ${\mathfrak p}_{-}( y_0^{T})$ and  ${\mathfrak p}_{+}( y_0^{T})$ respectively.  
Invoking again the path-wise control, $y_t^T \ge z_t$,  $t \in [-T, 0]$ we have that
$$
   {\mathfrak p}_1(y_0^T) e^{\psi(y^T)}   \le  {\mathfrak p}_{+}( z_0), \quad 
  {\bold{E}} [  {\mathfrak p}_{+}( z_0)] = \int_{-\infty}^0 {\mathfrak p}_{+}(z) \mathfrak{m}(dz) < \infty,
$$ 
where $\mathfrak{m}$ is defined in (\ref{zlaw}).
Returning to (\ref{elambda}), (\ref{pathwise}) and
 dominated convergence now produce
$$
  \lim_{\lambda \rightarrow \infty}  {\mathfrak e}_{\lambda} = \lim_{T \rightarrow \infty} {\bold{E}} [ {\mathfrak p}_1(y_0^T) e^{\psi(y^T)}   ]   
  = {\bold{E}} [  {\mathfrak p}_1(y_0^{\infty})  e^{\psi_{\infty} (y^{\infty})} ]  := \mathfrak{e},
$$
defining the constant $\mathfrak{e}$ in the statement of Theorem 1.  That $\mathfrak{e}$ is not equal to zero follows from
$$
\mathfrak{e} \ge \liminf_{T \rightarrow \infty} 
  {\bold{E}} [ {\mathfrak p}_1(y_0^T) e^{\psi(y^T)}   ] \ge \int_{-\infty}^0 {\mathfrak p}_{-}(z) \mathfrak{m}(dz) > 0. 
$$
Here we have used that $z \mapsto {\mathfrak p}_{-}(z)$ is decreasing in order to replace $y^T$-paths with $z$-paths, along with the fact that
${\mathfrak p}_1(z)$ (and so too  ${\mathfrak p}_{-}(z)$) is strictly positive (Lemma \ref{lem2}). 
 This completes the proof of Theorem 1,
granted the below.

\medskip

{\em Proof of Lemma 4.}
The continuity follows from that of the transition density $p(\cdot, x, y)$ in both space variables (the corresponding generator is hypo-elliptic).

\medskip

To see that ${\mathfrak p}_{1}(x) > 0$, first note that the operator $L_{\beta, a}$ defined in (\ref{L_integral})
which encodes the point process of eigenvalues  is positive and compact.  A proof that $L_{\beta, a}$ is in fact (almost surely) trace class 
is contained in Lemma 6 of \cite{RR}.  Its maximal eigenvalue, $\Lambda^{-1}$,  is therefore almost surely bounded above, and so there exists a small enough $\lambda_0 > 0$ such that ${\mathfrak p}_{\lambda_0}= {\mathfrak p}_{\lambda_0}(\infty)  > 0$.  Next, by the Markov property,
$$
                   {\mathfrak p}_{\lambda_0}(\infty)     = \int_{\-\infty}^{\infty} p(t, \infty, x) {\mathfrak p}_{\lambda_0 e^{-\beta t/4}} (x) dx, 
$$
and it follows that for every $t>0$ there is a $x_0$ such that $ {\mathfrak p}_{\lambda_0 e^{-\beta t/4}} (x_0) >0$.  
Applying the same formula
once again, we find that for any $z \in \RR$
$$
  {\mathfrak p}_{1}(x) \ge  \int_z^{\infty}   p(s,  x, y) {\mathfrak p}_{e^{-\beta s/4}} (y) dy    \ge  \PP_x( x_s \ge z)   {\mathfrak p}_{e^{-\beta s/4}}(z).
$$
To finish, choose $s= t -  \frac{4}{\beta} \log \lambda_0$ and then set $z$ to be the appropriate $x_0$. 

\medskip

For the bound (\ref{hittingbound}) we can restrict to $x$ less than some large negative constant, 
and note that ${\mathfrak p}_1(x)$ is bounded by the probability of non explosion for the following process
\begin{equation*}
d\tilde{y}_t = db_t + \frac{\beta}{4}\left( a+\frac{1}{2}\right) -  \frac{\beta}{4} e^{-\beta t/8} e^{-\tilde{y}_t}.
\label{eq:AnotherProcess}
\end{equation*}
since the downward drift on $\tilde{y}$ is weaker than that of $x$.
Next make the change $y_t = \tilde{y}_t + \frac{\beta t}{8}$ to obtain the homogenous process
\begin{equation*}
dy_t = db_t + \frac{\beta}{4}\left( a+1 - e^{-y_t}\right),
\label{eq:HomogProcess}
\end{equation*}
to which we can apply a version of Feller's test, similar to what was done at the end of Section 2.
A scale function for the $y$-process is 
\begin{equation*}
s(y) = \int_0^y \exp\{-\frac{\beta}{2} \left[ (a+1) \xi + e^{-\xi} -1 \right] \}\, d\xi,
\label{eq:Scale*}
\end{equation*}
and the probability that this process exits through $+\infty$ is exactly the probability of not exploding. 
This follows by checking the conditions 
of now Proposition 5.22 of \cite{KS}. According to that same proposition, the exit probability equals
\begin{equation*}
\frac{s(x)-s(-\infty)}{s(+\infty)-s(-\infty)} = \frac{1}{Z} \int_{-\infty}^x \exp\{-\frac{\beta}{2} \left[ (a+1) \xi + e^{-\xi}\right] \}\, d\xi,
\label{eq:ProbExit}
\end{equation*}
from which the required bound easily follows.
$\hfill$ $\square$

\bigskip

\noindent{\bf{Acknowledgements}}  The second and third named authors were supported in part by NSF grants DMS-0645756 and
 DMS-0804133, respectively.

\end{document}